\documentclass[a4paper]{article}
\pdfoutput=1
\usepackage[utf8]{inputenc}
\usepackage{amssymb}
\usepackage{biblatex}
\usepackage{amsmath}
\title{Non-trivial Poincare Constant in any Open, Bounded Subset $\Omega\subset \mathbb{R}^n$}
\author{Elliott Zaresky-Williams}
\date{November 18th, 2018}

\begin{document}

\maketitle

\begin{abstract}
    The Poincare Inequality is an extremely useful tool in the analysis of PDEs. A significant amount of literature has dealt with finding the optimal constant $C(p,\Omega)$, depending only the domain $\Omega$, and the $L^p$ norm in question. For convex subsets $\Omega\subset \mathbb{R}^n$, an optimal constant is known. There has also been some work done for \textit{specific} non-convex subsets. For a general open, bounded subset, finding the optimal constant is a non-trivial task. This paper will outline a method for finding a non-trivial estimate for $C$ when $\Omega\subset \mathbb{R}^n$ is non-convex, by first observing $\Omega\subset \mathbb{R}^n$ has finite Lebesgue measure, and using local estimates on subsets of $\Omega$ to obtain a global estimate on $\Omega$. 
\end{abstract}
\section{Introduction}

The classical statement of the Poincare inequality is as follows:

Let $\Omega\subset \mathbb{R}^n$ be an open, bounded, and connected subset and let $u\in W_0^{1,p}(\Omega)$. There exists $C(p,\Omega)$ where $||u||_{L^p(\Omega)}\leq C||\nabla u||_{L^p(\Omega)}$. In the case where $\Omega$ is convex, we can find an optimal $C(p,\Omega)$ for any $p$. In the case where $p=1$, the optimal constant is $C=diam(\Omega)/2$, which was shown in [1]. 

In the case for $p\geq 2$, the optimal constant $C(p,\Omega)$ is given by
$(\frac{\pi_p}{d})^p$, where $\pi_p=\frac{2\pi(p-1)^{1/p}}{psin(\pi/p)}$, which was shown in [2]. 

If we relax the assumption that $\Omega$ is convex, we still have an optimal estimate for $p=2$. 

The optimal constant $C(p,\Omega)$ is $\lambda^{-1}$ where $\lambda$ is the first Dirichlet eigenvalue, as this eigenvalue minimizes the Dirichlet energy $\frac{||\nabla u||_{L^2(\Omega)}}{||u||_{L^2(\Omega)}}$. 

\section{Familiar Estimates}

Since $\Omega\subset \mathbb{R}^n$ is bounded by hypothesis, we have that the Lebesgue measure $\mu(\Omega)<\infty$. Recall the standard estimate: 
\begin{center} $||u||_{L^p(\Omega)}\leq \mu(\Omega)^{1/p-1/q}||u||_{L^q(\Omega)}$ \end{center}

We can immediately relate this estimate to the original Poincare inequality. For sake of convenience, let $M=\mu(\Omega)^{1/p-1/q}$. Set $C(p,\Omega)=\frac{M||u||_{L^q(\Omega)}}{||\nabla u||_{L^p(\Omega)}}$. Notice that if we combine the inequalities, we have the following result: 

\begin{center} $||u||_{L^p(\Omega)}\leq C(p,\Omega)||\nabla u||_{L^p(\Omega)}\leq M||u||_{L^q(\Omega)}$ \end{center} 

\begin{center} (1) $||u||_{L^p(\Omega)}\leq \frac{M||u||_{L^q(\Omega)}}{||\nabla u||_{L^p(\Omega)}}||\nabla u||_{L^p(\Omega)}= M||u||_{L^q(\Omega)}$ \end{center} 

We now have a non-trivial estimate for any open, bounded, not necessarily convex subset $\Omega$. 

\section{Main Results}

The goal of the remainder of this section will be to use $(1)$, and the topology of $\Omega$ to find a non-trivial Poincare constant. First, local estimates will be given, and then the local estimates will be extended to the global result. 

The cases where $\Omega$ is \textit{not} bounded, and when $\Omega$ is \textit{totally} disconnected will also be discussed. 

How we partition $\Omega$ requires care. Any open subset of $\mathbb{R}^n$, can be written as a countable union of open balls $\Omega=\bigcup_{k}U_k$, but this is not the best partition of the set. If $\Omega$ is connected, the open balls cannot all be mutually disjoint. 

The intersection of each open ball contains another open ball which has positive measure. This would require splitting the theorem into two cases: one where $\Omega$ is connected and one where it is not. A "better" partition will be used, and the main theorem will address both $\Omega$ connected and disconnected simultaneously. 

Let $D_{j,k}=\bigcup_{j,k} [\frac{2^j}{k}, \frac{2^j}{k+1}]^n$ be the union of dyadic cubes in $\Omega$. Note that the interior $int(D_{j,k})=\bigcup_{j,k} (\frac{2^j}{k}, \frac{2^j}{k+1})^n$ has the same Lebesgue measure as $D_{j,k}$. We can apply $(1)$ to every open set $(\frac{2^j}{k}, \frac{2^j}{k+1})^n$ in the interior.

\subsection{Lemma 1}

Denote the interior of $D_{j,k}$ by $I_{j,k}$. Let $u\in L^p(D_{j,k})$ with $||u||_{L^p(D_{j,k})}=M$, then $u\in L^p(I_{j,k})$ with $||u||_{L^p(D_{j,k})}=||u||_{L^p(I_{j,k})}=M$.

$\forall j,k$ we have either than $D_{j,k}\subset D_{j+1,k+1}$  or $\mu(D_{j,k}\bigcap D_{j+1,k+1})=\emptyset$. It immediately follows that $\mu(D_{j,k})=\mu(I_{j,k})$ and that  $$||u||_{L^p(D_{j,k})}=||u||_{L^p(I_{j,k})}$$

There was nothing particularly special about choosing dyadic cubes except for familiarity and simplicity. In general, we can replace the open dyadic cubes with any collection of open sets whose almost disjoint union is $\Omega$. Here, and the remainder of the paper, two sets $U,V$ are "almost disjoint" if their intersection has measure zero (i.e. $\mu(U\bigcap V)=0$).

\subsection{Lemma 2}

Let $\Omega=\bigcup_k U_k$ be a countable, almost disjoint union of open sets. If $u\in L^p(\Omega)$, then we have $||u||^p_{L^p(\Omega)}=\sum_k ||u||^p_{L^p(U_k)}$.

This result builds off the countable subadditivity of the Lebesgue integral for $L^1$ functions. Namely if $u\in L^1(\Omega)$ then we have
$\int_{\Omega}u d\mu=\sum_k \int_{U_k}ud\mu$ where $d\mu$ denotes the Lebesgue measure. 

The proof of the above is standard. For the lemma, make the substitution $g=|u|^p$, and we arrive at $$\int_{\Omega}g d\mu=\sum_k \int_{U_k} g d\mu$$ Now, raise everything to $1/p$ to obtain $$(\int_{\Omega}g d\mu)^{1/p}=(\sum_k \int_{U_k} g d\mu)^{1/p}$$
Now we can back substitute $|u|^p$, arriving at $$(\int_{\Omega}|u|^p d\mu)^{1/p}=(\sum_k \int_{U_k} |u|^p d\mu)^{1/p}$$ Finally, raising everything to the power of $p$ we obtain $$(\int_{\Omega}|u|^p d\mu)=(\sum_k \int_{U_k} |u|^p d\mu)$$ which is just  $$||u||^p_{L^p(\Omega)}=\sum_k ||u||^p_{L^p(U_k)}$$ and the lemma follows. 

\subsection{Lemma 3}

Let $\Omega=\bigcup_k U_k$ be a countable (almost disjoint) union of bounded open sets, and assume $p\geq 1$. We have that $u\in L^p(\Omega)$ if and only if $u\in L^p(U_k)$ $\forall k$ and $\sum_k ||u||^p_{L^p(U_k)}\leq \sum_k \frac{1}{k^{1+\epsilon}}$ for some $\epsilon>0$. 

$\implies)$ Since $u\in L^p(\Omega)$, we have $||u||_{L^p(\Omega)}=N$ (since norm is finite). By lemma 2, we have $$||u||^p_{L^p(\Omega)}=\sum_k ||u||^p_{L^p(U_k)}$$ Continuing, $$||u||^p_{L^p(\Omega)}=N^p=\sum_k ||u||^p_{L^p(U_k)}$$ 

$$N^p=\sum_k ||u||^p_{L^p(U_k)}\leq \sum_k \frac{1}{k^{1+\epsilon}}$$ 

The series $\sum_k \frac{1}{k^{1+\epsilon}}$ is convergent $\forall \epsilon>0$, so we just need to find an $\epsilon>0$ which satisfies the above. 

The series $\sum_k \frac{1}{k^{1+\epsilon}}$ is known to be upper bounded by 
$K_{\epsilon}=\frac{1+\epsilon}{\epsilon}$. We now just need 

$$N^p\leq \frac{1+\epsilon}{\epsilon}$$

We now just pick $\epsilon\leq \frac{1}{N^p-1}$, and we are done. 

$\impliedby)$ For the converse, we combine lemma 2 and the hypothesis to arrive at $$||u||^p_{L^p(\Omega)}=\sum_k ||u||^p_{L^p(U_k)}\leq \sum_k \frac{1}{k^{1+\epsilon}}$$ for some $\epsilon>0$. 

In particular, we have $$||u||^p_{L^p(\Omega)}\leq \sum_k \frac{1}{k^{1+\epsilon}}$$ After raising to the $1/p$ power, we obtain

$$||u||_{L^p(\Omega)}\leq (\sum_k \frac{1}{k^{1+\epsilon}})^{1/p}\leq \sum_k \frac{1}{k^{1+\epsilon}}$$

Recall the last inequality is justified since $p\geq 1$ by hypothesis. 

The condition of $\sum_k ||u||^p_{L^p(U_k)}\leq \sum_k \frac{1}{k^{1+\epsilon}}$ cannot be dispensed with, as the following counterexample shows. 

Let $\Omega=(0,1)$ with $U_k=(1/k,1)$, clearly $\Omega=\bigcup_k U_k$. Let $f=1/x \chi_{(1/k,1)}$.

We have  $$||f||_{L^1(U_k)}=-ln(1/k)$$ clearly $f\in L^1(U_k)$, but the sum $\sum_k ||f||_{L^1(U_k)}=-ln(1/k)$ is unbounded and $f\notin L^1(\Omega)$. In general, if the sum of the $L^p(U_k)$ norms form an unbounded series, the function $f$ will not lie in $L^p(\Omega)$.

\subsection{Theorem}

Let $\Omega\subset \mathbb{R}^n$ be an open, bounded subset where $\Omega=\bigcup_{k} U_k$ a countable union of open sets which are bounded, open, and pairwise almost disjoint. 

For $1\leq p<q\leq \infty$, let $u\in W^{1,p}(U_k)\bigcap L^q(U_k)$ $\forall k$. Finally, assume $\sum_k ||\nabla u||^p_{L^p(U_k)}\leq \sum_k \frac{1}{k^{1+\epsilon}}$ and $\sum_k \mu(U_k)^{1-p/q}||u||^p_{L^q(U_k)}\leq \sum_k \frac{1}{k^{1+\delta}}$ for some $\epsilon,\delta>0$, then a non-trivial Poincare constant $C(p,\Omega)$ can be found.

Note that the open dyadic cubes of the form $I_{j,k}=(\frac{2^j}{k}, \frac{2^j}{k+1})^n$ can be chosen, as they fit the criteria of the theorem, as the set $I_{j,k}$ is at most countable. 

By lemma 3, and the hypothesis that $u\in W^{1,p}(U_k)$ $\forall k$, we have $$||\nabla u||^p_{L^p(U_k)}\leq \sum_k \frac{1}{k^{1+\epsilon}}$$ and
$$||\nabla u||^p_{L^p(U_k)}=||\nabla u||^p_{L^p(\Omega)}$$

together, this implies
$$||\nabla u||_{L^p(\Omega)}\leq (\sum_k \frac{1}{k^{1+\epsilon}})^{1/p} \leq \sum_k \frac{1}{k^{1+\epsilon}}$$. 

Similarly, since $u\in L^q(U_k)$, $\forall k$, we have 
$$||u||_{L^q(\Omega)}\leq (\sum_k \frac{1}{k^{1+\delta}})^{1/p}\leq \sum_k \frac{1}{k^{1+\delta}}$$

Combing these two allows us to ensure that both $||\nabla u||_{L^p(\Omega)}$ and $||u||_{L^q(\Omega)}$ are bounded above. 

By lemma 2, we obtain:

$(i)  \hspace{4mm} 
||u||^p_{L^p(\Omega)}=\sum_k ||u||^p_{L^p(U_k)}$

$(ii)\hspace{3mm} ||\nabla u||^p_{L^p(\Omega)}=\sum_k||\nabla u||^p_{L^p(U_k)}$ 

We simply rewrite these quantities as 
 
$(i) \hspace{4mm} ||u||_{L^p(\Omega)}=(\sum_k ||u||^p_{L^p(U_k)})^{1/p}$

$(ii) \hspace{3mm} ||\nabla u||_{L^p(\Omega)}=(\sum_k ||\nabla u||^p_{L^p(U_k)})^{1/p}$

Care must be taken with extending the local $L^q$ estimate to the global one. Locally, for each $U_k$, we have the estimate 
$$||u||_{L^p(U_k)}\leq \mu(U_k)^{1/p-1/q}||u||_{L^q(U_k)}$$ 

However, raising the above quantity to the power of $p$, we obtain 

$$||u||^p_{L^p(U_k)}\leq \mu(U_k)^{1-p/q}||u||^p_{L^q(U_k)}$$ Summing up the local estimates, we have

$$||u||^p_{L^p(\Omega)}=\sum_k||u||^p_{L^p(U_k)}\leq \sum_k \mu(U_k)^{1-p/q}||u||^p_{L^q(U_k)}$$

The rightmost sum is convergent by hypothesis, so there is no issue extending the local $L^q$ estimate to a global one. 

We need an exact value for the global $L^q$ norm. The exact value follows from lemma 2:

$$||u||^q_{L^q(\Omega)}=\sum_k ||u||^q_{L^q(U_k)}$$ which becomes
$(iii) \hspace{10mm} ||u||_{L^q(\Omega)}=(\sum_k ||u||^q_{L^q(U_k)})^{1/q}$

The last quantity we need is to relate $\mu(\Omega)$ to $\mu(U_k)$. By countable subadditivity of the Lebesgue measure, we have that $\mu(\Omega)=\sum_k \mu(U_k)$, since $\Omega=\bigcup_k U_k$. We arrive at

$(iv) \hspace{2mm} \mu(\Omega)^{1/p-1/q}=(\sum_k \mu(U_k))^{1/p-1/q}$

What remains now is to tie these inequalities together into the statement of the Poincare inequality. 

The Poincare constant we desire is $$C(p,\Omega)=\frac{\mu(\Omega)^{1/p-1/q}||u||_{L^q(\Omega)}}{||\nabla u||_{L^p(\Omega)}}$$

After substituting the the sums of local estimates $(ii-iv)$ in place of global ones, we obtain 

$$C(p,\Omega)=\frac{(\sum_k \mu(U_k))^{1/p-1/q}(\sum_k ||u||^q_{L^q(U_k)})^{1/q}}{(\sum_k ||\nabla u||^p_{L^p(U_k)})^{1/p}}$$

\section{Conclusion} 

Here we see why we need to partition $\Omega$ in such a way. If the intersection of the union of the almost disjoint $U_k$ had non-zero Lebesgue measure, lemma 3 would not hold, and the lemma was vital to the proof. 

Note also that no assumption was made about the connectedness of $\Omega$. The theorem still holds, even if $\Omega$ was totally disconnected (but still with positive Lebesgue measure).

In addition, the weakening of the assumption that $u\in W^{1,p}_0(U_k)$ to just $W^{1,p}(U_k)$ is justified. The $bd(U_k)$ has measure zero. Even if $u$ doesn't vanish on the boundary, the norm is not affected. 

Unfortunately, when $\Omega\subset \mathbb{R}^n$ is \textit{unbounded}, we cannot obtain a \textit{general} estimate of the form 

$$||u||_{L^p(\Omega)}\leq C(p,\Omega)||\nabla u||_{L^p(\Omega)}$$

The reasoning is straightforward. Observe that for the bounded case, we have
$\mu(\Omega)=M$ (finite measure), and that the Poincare constant was 
$$C(p,\Omega)=\frac{(\sum_k \mu(U_k))^{1/p-1/q}(\sum_k ||u||^q_{L^q(U_k)})^{1/q}}{(\sum_k ||\nabla u||^p_{L^p(U_k)})^{1/p}}$$

The sum $(\sum_k \mu(U_k))^{1/p-1/q}$ is the problematic term. The whole space $\mathbb{R}^n$ is second countable, and so in particular $\Omega$ is second countable. We can find a countable topological base of open subsets whose union is $\Omega$. The issue is that it is not possible to find such a union of open sets where the sum $(\sum_k \mu(U_k))^{1/p-1/q}$ converges, because if it did, then $\Omega$ would have finite measure, contradicting that $\Omega$ was unbounded.

\end{document}